\documentclass[11pt,letterpaper]{amsart}

\usepackage[english]{babel}
\usepackage[OT2,T1]{fontenc} 

\usepackage{mathtools}
\usepackage{amssymb}
\usepackage{stmaryrd}
\usepackage{bm}
\usepackage{tikz-cd}

\usepackage{comment}
\usepackage{todonotes}
\usepackage{xcolor}
\usepackage{enumitem}
\usepackage{graphicx}
\usepackage{setspace}
\usepackage{url}
\usepackage{xr}

\usepackage[title]{appendix}
\usepackage{chngcntr}
\usepackage{apptools}
\usepackage{thmtools}

\usepackage[pagebackref]{hyperref}
\usepackage[nameinlink]{cleveref}
\usepackage{bookmark}

\definecolor{crimsonglory}{rgb}{0.75, 0.0, 0.2}
\definecolor{darkpowderblue}{rgb}{0.0, 0.2, 0.6}

\hypersetup{
	colorlinks = true,
	urlcolor   = darkpowderblue,
	linkcolor  = darkpowderblue,
	citecolor  = crimsonglory
}

\theoremstyle{plain}
\newtheorem{theorem}{Theorem}[section]
\newtheorem{definition}[theorem]{Definition}

\newtheorem{prop}[theorem]{Proposition}
\newtheorem{cor}[theorem]{Corollary}
\newtheorem{lemma}[theorem]{Lemma}
\newtheorem{claim}[theorem]{Claim}
\newtheorem{remark}[theorem]{Remark}

\newtheorem{conj}[theorem]{Conjecture}

\DeclareMathOperator{\id}{id}

\DeclareMathOperator{\End}{End}
\DeclareMathOperator{\GL}{GL}

\DeclareMathOperator{\tr}{tr}

\DeclareMathOperator{\spec}{Spec}
\DeclareMathOperator{\gal}{Gal}

\DeclareMathOperator{\Span}{Span}

\DeclareMathOperator{\disc}{disc}

\DeclareMathOperator{\crys}{crys}

\DeclareMathOperator{\ord}{ord}
\DeclareMathOperator{\SP}{Sp}
\DeclareMathOperator{\adj}{adj}
\DeclareMathOperator{\arch}{arch}

\DeclareMathOperator{\ssing}{ssing}

\DeclareMathOperator{\an}{an}
\DeclareMathOperator{\trdeg}{trdeg}

\newcommand{\Sum}[2]{\displaystyle\sum_{#1}^{#2}}

\newcommand{\N}{\mathbb{N}}
\renewcommand{\C}{\mathbb{C}}
\newcommand{\Q}{\mathbb{Q}}
\newcommand{\B}{\mathbb{B}}
\newcommand{\A}{\mathbb{A}}

\newcommand{\CX}{\mathcal{X}}

\newcommand{\CP}{\mathcal{P}}
\newcommand{\CO}{\mathcal{O}}

\newcommand{\CF}{\mathcal{F}}

\DeclareSymbolFont{cyrletters}{OT2}{wncyr}{m}{n}
\DeclareMathSymbol{\Sha}{\mathalpha}{cyrletters}{"58}

\newcommand{\mathsym}[1]{{}}
\newcommand{\unicode}[1]{{}}

\title{$G$-functions and quaternionic multiplication in $\mathcal{A}_2$}
\author{Georgios Papas}
\address{Faculty of Mathematics and Computer Science\\
	The Weizmann Institute of Science\\
	234 Herzl Street\\
	Rehovot 76100\\
	Israel, and\\
	Institute for Advanced Study\\
	1 Einstein Drive\\
	Princeton, N.J. 08540\\
	U.S.A.}

\email{georgios.papas@weizmann.ac.il}
\email{gpapas@ias.edu}

\begin{document}
	
\begin{abstract}
	We consider $G$-functions associated with one--parameter families of principally polarized abelian surfaces. We construct relations among their values at QM points, under the assumption that the family is Hodge generic and has a distinguished QM point.
	
	The proof uses archimedean and $p$-adic period matrices of QM abelian surfaces, via the de Rham--Betti and de Rham--crystalline comparison isomorphisms. As an application of our relations and the Andr\'e--Bombieri Hasse principle, we obtain a height bound for QM points and discuss consequences for the Zilber--Pink conjecture in $\mathcal{A}_2$.
\end{abstract}
	
	\maketitle
	\section{Introduction}\label{section:intro}

We study values of $G$-functions attached to one--parameter families of principally polarized abelian surfaces. We consider families with a distinguished point whose fiber has quaternionic multiplication and study the values of the associated $G$-functions at points whose fibers also have quaternionic multiplication.

More concretely, let $f:\CX\rightarrow S$ be a one--parameter family of principally polarized abelian surfaces defined over a number field $K$, and fix $s_0\in S(K)$ such that the fiber $\CX_{s_0}$ has quaternionic multiplication\footnote{Throughout, we say that an abelian surface has quaternionic multiplication, or is QM for short, if its geometric endomorphism algebra is an indefinite quaternion algebra over $\Q$.}. Associated with this data is a family $Y_G$ of $G$-functions, constructed from a basis of sections of $H^1_{dR}(\CX/S)$ over an affine open neighborhood $U$ of $s_0$; see \Cref{section:backgroundgfuns}.

For an algebraic point $s\in S(\overline{\Q})$, we say that a place $v$ of $K(s)$ is a place of proximity of $s$ to $s_0$ if it satisfies the conditions of \Cref{section:backgroundgfuns}. At every such place, the corresponding $G$-functions converge at $s$, and we write $Y_G(s)$ for their values in $\C_v$.

For a place $v$ of a number field $L$, we fix an embedding $\iota_v:L\hookrightarrow\C_v$ inducing $v$. Our main result is the following.

\begin{theorem}\label{proprelationsintro}
	Assume that the induced morphism $\iota:S\rightarrow \mathcal{A}_2$ has image a Hodge generic curve. 
	
	Let $s\in S(\overline{\Q})$ be a point such that $\CX_s$ has quaternionic multiplication. There exists a number field $L_s\supseteq K(s)$, with $[L_s:K(s)]$ bounded independently of $s$, such that, for every place $v$ of $L_s$ above a place of proximity of $s$ to $s_0$, there is a polynomial $R_{s,v}\in L_s[\underline X]$ of degree at most $8$ satisfying
	$$
	\iota_v\bigl(R_{s,v}(Y_G(s))\bigr)=0.
	$$
	
	This relation is non-trivial, i.e. it does not arise by specializing at $s$ a functional relation among the entries of $Y_G$.
	
	Moreover, at places $v$ of $L_s$ lying above places of proximity at which $\CX_{s_0}$ has potentially ordinary reduction, the polynomial $R_{s,v}$ may be chosen independently of $v$.
\end{theorem}

When the point $s_0$ corresponds to an interior point of the relevant moduli space rather than a boundary point, algebraic relations among values of $G$-functions, either archimedean or $p$-adic, are much less understood. For one--parameter families of elliptic curves, such relations at CM points were studied by Andr\'e \cite{andremots,andremsj} and Beukers \cite{beukers}. Andr\'e's approach is cohomological, whereas Beukers obtains the $p$-adic relations by lifting isogenies of the reduction to isogenies of families in characteristic zero. The local picture for such relations was completed in \cite{papaspadicpart2}, where relations were constructed at every place of proximity.

Following Andr\'e's cohomological approach \cite{andremots,andremsj}, we use the de Rham--Betti and de Rham--crystalline comparison isomorphisms to interpret the values of $Y_G$ in terms of period matrices of the fibers of $f$. Our cohomological framework applies uniformly at all places of proximity, although the period--theoretic input differs in the ordinary and supersingular cases.

Applications of the $G$-functions method to unlikely intersections have largely concerned distinguished points $s_0$ at which the family degenerates, corresponding to intersections with the boundary of the ambient moduli space; see \cite{daworr4,daworr5,papasbigboi,papasimrn,urbanik}. In that setting, the geometry of the degeneration yields $p$-adic relations that are uniform in the place of proximity. When $s_0$ corresponds to an interior point, no analogous global geometric source of relations is available, and arithmetic information at all finite places becomes essential.

\subsection{Outline of the relation construction}

At every finite place, a QM abelian surface has potentially good reduction and is either potentially ordinary or potentially supersingular; see \Cref{section: qm and reductions}. This dichotomy is closely analogous to that for CM elliptic curves. Together with the archimedean case, it gives the three cases that organize our construction.

\subsubsection{Period matrices of QM abelian surfaces}

We obtain explicit forms for the period matrices of QM abelian surfaces in both the archimedean and non--archimedean settings by transporting the quaternionic action through the corresponding comparison isomorphisms.

In the archimedean case, we choose an imaginary quadratic subfield of the quaternion algebra on which the Rosati involution acts by complex conjugation; see \Cref{qmalglemma1}. This induces eigenspace decompositions on both sides of the de Rham--Betti comparison isomorphism, and choosing bases adapted to these decompositions yields the required form of the period matrix; see \Cref{qmarch}.

The same imaginary quadratic subfield, and the corresponding decomposition of de Rham cohomology, are used at finite places. On the crystalline side, the two reduction types require different treatment. 

At potentially ordinary places, compatibility of the endomorphism action with Frobenius yields greater rigidity: the Frobenius slope decomposition is preserved by the endomorphism action, and its interaction with the eigenspace decomposition forces a more restrictive form of the period matrix; see \Cref{qmnonarch}. Exploiting this slope decomposition is a key new ingredient in our construction. This additional rigidity allows the relation in \Cref{proprelationsintro} to be chosen independently of the place in the potentially ordinary case.

In the potentially supersingular case, Frobenius no longer has distinct slopes. Andr\'e's $p$-adic Betti lattices therefore replace the Frobenius slope decomposition as the principal tool on the crystalline side; see \Cref{qmsupersinglembasis} and \Cref{qmsupersingperiod0}.

\subsubsection{From periods to values of $G$-functions}

The cohomological interpretation of the power series $Y_G$ is provided by the relative de Rham--Betti and de Rham--crystalline comparison isomorphisms. For algebraic points $s$ lying in a sufficiently small analytic disk around $s_0$, in the complex or rigid analytic sense, the period matrix of $\CX_s$ can be written as
$$
\Pi_v(\CX_s)=\iota_v\bigl(Y_G(s)\bigr)\cdot\Pi_v(\CX_{s_0}),
$$
where $\Pi_v(\CX_{s_0})$ denotes the period matrix of the distinguished fiber; see \Cref{section:backgroundgfuns}.

Thus, algebraic relations among the entries of $\Pi_v(\CX_s)$ translate into algebraic relations among the values $Y_G(s)$, after eliminating the transcendental periods of $\CX_{s_0}$. The explicit descriptions of the period matrices obtained above allow us to carry out this elimination.

\subsubsection{Non-triviality of the relations}

It remains to show that these relations are non-trivial, i.e. that they are not specializations at $s$ of functional relations satisfied by the entries of $Y_G$.

Under the Hodge genericity assumption of \Cref{proprelationsintro}, work of Bakker--Tsimerman \cite{bakkertsimerman} gives a complete description of these functional relations. With the conventions used in \Cref{section:backgroundgfuns}, the ideal of functional relations is precisely $I(\SP_4)$, the ideal defining $\SP_4$ as a subvariety of $\GL_4$; see \cite{papaspadicpart1}.

The accompanying Mathematica code \cite{mathematicacodeqm} computes the relevant coefficient expansions of the polynomials used in the proof of \Cref{proprelationsintro} and their remainders upon division by a Gr\"obner basis of $I(\SP_4)$. We use these computations in \Cref{section:qmfiniteplaces} to establish the required non-vanishing statements, and hence the non-triviality of the constructed relations.

\subsection{Height bounds and applications}

As an application of the relations constructed in \Cref{proprelationsintro}, the Andr\'e--Bombieri Hasse principle for values of $G$-functions yields a height bound for points whose fibers have quaternionic multiplication.

For $s\in S(\overline{\Q})$, let $\Sigma_{\ssing}(s,s_0)$ be the set of places $w\in\Sigma_{K,f}$ at which $\CX_{s_0}$ has potentially supersingular reduction and for which there exists $v\in\Sigma_{K(s),f}$ with $v\mid w$ such that $s$ is $v$-adically close to $s_0$.
\begin{theorem}\label{Theoremmainqmhtbound}	
	Let $f:\CX\rightarrow S$ be a one--parameter family of principally polarized abelian surfaces as in \Cref{proprelationsintro}, and let $h$ be a Weil height on $S$. 
	
	There exist positive constants $c_1$, $c_2$, depending only on the family and on $s_0$, such that
	$$
	h(s)\leq c_1\Bigl([\Q(s):\Q]\bigl(1+|\Sigma_{\ssing}(s,s_0)|\bigr)\Bigr)^{c_2},
	$$
	for all points $s\in S(\bar{\Q})$ for which $\CX_s$ has quaternionic multiplication. 
\end{theorem}

The appearance of the term $|\Sigma_{\ssing}(s,s_0)|$ in the height bound reflects the dependence of $R_{s,v}$ on the place $v$ in \Cref{proprelationsintro}. No corresponding term arises from potentially ordinary places, since there the relation may be chosen independently of $v$. We return to this point in \Cref{section:placesofssingprox}. 

Our height bound may be viewed as a natural analogue of the statement announced as Theorem~1 in the introduction to \cite{andremots}. In the modular elliptic setting, Andr\'e's statement concerns CM points and is thus closer in spirit to Andr\'e--Oort, whereas our result concerns QM points on a Hodge generic curve in $\mathcal{A}_2$. We discuss the connection of our height bound with the Zilber--Pink conjecture in $\mathcal{A}_2$ in \Cref{section: zilberpinkapplications}.

\subsection{Notation}
Let $X$ be an abelian variety over a number field $K$ and $v\in\Sigma_{K}$ a place of $K$. We write $X_v$ for the base change $X\times_{K} K_v$. If $v$ is a finite place of good reduction for $X$, we write $\tilde{X}_v$ for the reduction of $X$ modulo $v$.

Given a family of power series $\mathcal{Y}:=(y_1,\ldots,y_N)\in K[[x]]$ with coefficients in a number field $K$, for any place $v\in\Sigma_K$ we write $R_v(y_j)$ for the $v$-adic radius of convergence of $y_j$. Similarly, for the $v$-adic radius of the entire family, we write $R_{v}(\mathcal{Y}):=\min R_v(y_j)$. For each such $v$, we fix an embedding $\iota_v:K\hookrightarrow\C_v$, where $\C_v$ denotes either $\C$ or $\C_p$ according as $v\mid\infty$ or $v\mid p$. In this spirit, for a power series $y(x)=\Sum{n=0}{\infty}a_n x^n\in K[[x]]$ we write $\iota_v(y(x)):=\Sum{n=0}{\infty}\iota_v(a_n)x^n\in\C_v[[x]]$.

For a finite place $v$ of residue characteristic $p$, we write $k_v$ for the residue field and $K_{v,0}:=W(k_v)[1/p]$.

	\section{QM abelian surfaces}\label{section:backqm}

We recall the facts about QM abelian surfaces used in the sequel, together with background on the quaternion algebras arising as their endomorphism algebras.

Throughout this section, all QM abelian surfaces are assumed to be equipped with a principal polarization defined over the ground field.

\subsection{Endomorphisms and reductions}\label{section: qm and reductions}

All quaternion algebras considered below are indefinite over $\Q$. Following \cite{voight}, we write
$$
D=\left(\frac{a,b}{\mathbb{Q}}\right),
$$
where $a,b\in\mathbb{Q}$ and $D$ has basis $\{1,i,j,k=ij\}$ satisfying
$$
i^2=a,\qquad j^2=b,\qquad k=ij=-ji.
$$
Whenever we use this notation, we assume that $a$ and $b$ are nonsquare integers; see \cite[Exercise~2.4]{voight}.

\begin{lemma}\label{qmalglemma1}
	Let $A$ be a principally polarized abelian surface over a number field $K$ and assume that $D:=\End^{0}_{K}(A)=\End^{0}_{\bar{\Q}}(A)$ is an indefinite quaternion algebra with $D=\left(\frac{a, b}{\mathbb{Q}}\right)$ and basis $\{1, i, j,k= ij\}$. Then there exists a pure quaternion\footnote{A pure quaternion in $D$ is a quaternion that can be written as $Xi+Yj+Zk$ with $X$, $Y$, $Z\in\Q$.} $\lambda\in D$ with $\lambda^2<0$ such that 
	$$
	\lambda^{\dagger}=-\lambda
	$$
	for the Rosati involution $\dagger$ on $D$.
	\end{lemma}
	\begin{proof}
		By \cite[Theorem~2, p.~186]{mumfordabelian}, there exists a pure quaternion $\alpha\in D$ such that $\alpha^2\in\Q$, $\alpha^2<0$, and the Rosati involution on $D$ is given by
		$$
		x^\dagger=\alpha x^*\alpha^{-1},
		$$
		where $x^*$ denotes the standard involution of $D$. Since $\alpha$ is pure, we have $\alpha^*=-\alpha$. Hence
		$$
		\alpha^\dagger=\alpha\alpha^*\alpha^{-1}=-\alpha.
		$$
		Taking $\lambda=\alpha$ proves the claim.
	\end{proof}
We next recall a theorem of Patankar describing the reductions of QM abelian surfaces.
\begin{theorem}[{\cite[Theorem~2.2.2]{patankarthesis}}]\label{reductionqmsurface}
		Let $A$ be an abelian surface defined over a number field $K$. Assume that $\End_K^{0}(A)=\End^{0}_{\bar{\Q}}(A)$ is a quaternion algebra and that $A$ has everywhere good reduction. Then for each finite place $v\in\Sigma_{K,f}$ there exists an elliptic curve $\tilde{E}_v$ defined over $k_v$, the residue field of $K$ at $v$, such that 
		\begin{center}
			$\tilde{A}_v\sim_{k_v} \tilde{E}_v^2$.
		\end{center}
	\end{theorem}
	
\begin{definition}
	A finite place $v$ is called \textbf{ordinary for} $A$ (respectively \textbf{supersingular for} $A$) if the associated elliptic curve $\widetilde{E}_v$ is ordinary (respectively supersingular).
\end{definition}

\subsubsection{Choice of the ground field}\label{section:conventionsqm}
	
	We next record a series of properties of QM abelian surfaces required in the sequel, which may be arranged simultaneously after a finite extension of the ground field.
	\begin{lemma}\label{qmlem1}
		Let $A$ be a principally polarized QM abelian surface defined over some number field $K$ such that $\End^{0}_{\bar{\Q}}(A)=:D$. Then there exists a finite extension $L/K$ such that the base change $A_L:=A\times_{\spec(K)} \spec(L)$ satisfies the following properties:
		\begin{enumerate}
			\item\label{qmlem1item1} $\End^{0}_L(A_L)=D$, 
			
			\item\label{qmlem1item2} for a choice of $\lambda$ as in \Cref{qmalglemma1}, we have $\Q(\sqrt{\lambda^2})\subset L$,
			
			\item\label{qmlem1item3} $A_L$ has everywhere good reduction.\end{enumerate}
	\end{lemma}
	\begin{proof}
		It is classical that, after replacing $K$ by a finite extension $L$, one has $\End_L^0(A_L)=D$; moreover, the degree $[L:K]$ of the extension may be bounded independently of $A$. See, for example, \cite{silverberg}.
		
		After replacing $L$ by its compositum with $\Q(\sqrt{\lambda^2})$, property~\ref{qmlem1item2} holds. Finally, the claim in the proof of \cite[Theorem~1.1]{clarkxarles} shows that $A$ has potentially good reduction at every finite place. After replacing $L$ by a further finite extension, Grothendieck's semistable reduction theorem \cite[Th\'eor\`eme~3.6]{sga} shows that $A$ has everywhere semistable reduction, which is necessarily good in view of potential good reduction. This proves \ref{qmlem1item3}.
	\end{proof}
	
	\subsection{Period matrices of QM abelian surfaces}\label{section:qmperiods}
	
	Let $A$ be a QM abelian surface over a number field $K$ satisfying the conclusions of \Cref{qmlem1}. In particular, $A$ has everywhere good reduction. For each $v\in\Sigma_K$, set
	\begin{equation}
		H^1_v(A)=	\begin{cases}	H^1(A_v^{an},\Q) & \text{if } v\in\Sigma_{K,\infty},\\
			H^1_{\crys}(\widetilde{A}_v/W(k_v))\otimes K_{v,0}
			& \text{if } v\in\Sigma_{K,f}.	\end{cases}
	\end{equation}
	We denote by $\rho_v(A)$ the corresponding comparison isomorphism: the de Rham--Betti comparison isomorphism when $v$ is archimedean, and the de Rham--crystalline comparison isomorphism when $v$ is finite; see \cite[Theorem~2.4]{bertogus}. Given suitable bases of $H^1_{dR}(A/K)$ and $H^1_v(A)$, we call the matrix representing $\rho_v(A)$ the period matrix of $A$ at $v$.
	
	We use the following class of bases of $H^1_{dR}(A/K)$ throughout the sequel.
	\begin{definition}\label{hodgebasis}Let $A$ be a principally polarized abelian surface over a number field $K$.  We call an ordered basis $\Gamma_{dR}(A):=(\omega_1,\omega_2,\eta_1,\eta_2)$ of $H^1_{dR}(A/K)$ a \textbf{Hodge basis} if the following conditions hold:
		\begin{enumerate}
			\item $\omega_1,\omega_2$ form a basis of the first step $F^1=e^*\Omega_{A/K}\subset H^1_{dR}(A/K)$ of the Hodge filtration, and
			
			\item $\Gamma_{dR}(A)$ is symplectic with respect to the de Rham polarization pairing, meaning that $\langle\omega_i,\eta_j\rangle=\delta_{i,j}$ and $\langle\omega_i,\omega_j\rangle=\langle\eta_i,\eta_j\rangle=0$.
		\end{enumerate}
	\end{definition}
	
	For the remainder of this section we adopt the notation and conclusions of \Cref{qmlem1}. In particular, we assume that $D:=\End^{0}_K(A)=\End^{0}_{\bar{\Q}}(A)$
	is an indefinite quaternion algebra. Fix a presentation
	$$
	D=\left(\frac{a,b}{\Q}\right)
	$$
	and a standard basis $\{1,i,j,ij\}$. Choose $\lambda$ as in \Cref{qmalglemma1}, and set $\CF_\lambda:=\Q(\lambda)\subset D$. We identify $\CF_\lambda$ with the imaginary quadratic field $\Q(\sqrt{\lambda^2})$ via $\lambda\mapsto\sqrt{\lambda^2}$. By our choice of the ground field in \Cref{qmlem1}, $\CF_\lambda\subset K$.
	
	We repeatedly use the following class of Hodge bases.
	\begin{definition}
		Let $F\subset\End^0(A)$ be a quadratic field with $F\subset K$, and write
		$$
		\gal(F/\Q)=\{\sigma,\tau\},
		$$
		where $\sigma=\id_F$ and $\tau$ is the nontrivial automorphism of $F/\Q$.
		We call an ordered Hodge basis $(\omega,\omega',\eta,\eta')$ of $H^1_{dR}(A/K)$ an \textbf{$F$-Hodge basis} if, for every $a\in F$,
		\begin{enumerate}
			\item $a^*\omega=\sigma(a)\omega$ and
			$a^*\eta'=\sigma(a)\eta'$, and
			\item $a^*\omega'=\tau(a)\omega'$ and
			$a^*\eta=\tau(a)\eta$.
		\end{enumerate}
		Here $a^*$ denotes the pullback action induced by $a$ on $H^1_{dR}(A/K)$.
	\end{definition}
	
	We now describe the period matrices of $A$ separately at archimedean, ordinary, and supersingular places. In each case, we choose bases of $H^1_v(A)$ adapted to the action of the endomorphism algebra, so that the corresponding period matrices admit explicit block descriptions.
	
	In the crystalline setting, these choices also give precise control over the matrices representing the pullback action of endomorphisms: at ordinary places these matrices admit an explicit description, while at supersingular places their coefficients lie in $\bar{\Q}$ rather than merely in $\C_p$.
	
	\subsubsection{Archimedean places}\label{section:qmperiodsarch}
	
	Let $v\in\Sigma_K$, let $\Gamma_v(A)$ be a symplectic basis of $H^1_v(A)$, and let $\Gamma_{dR}(A)$ be a Hodge basis of $H^1_{dR}(A/K)$. We write
	\begin{equation}
		\Pi_v(\Gamma_{dR}(A),\Gamma_v(A)):=[\rho_v(A)]_{\Gamma_v(A)}^{\Gamma_{dR}(A)}
	\end{equation}
	for the period matrix of $A$ with respect to these bases.
	
	\begin{lemma}\label{qmarch} 
		Let $v\in\Sigma_{K,\infty}$. Then one can choose an $\CF_{\lambda}$-Hodge basis	$\Gamma_{dR}(A):=(\omega,\omega',\eta,\eta')$ of $H^1_{dR}(A/K)$ and a symplectic basis $\Gamma_v(A)$ of $H^1_v(A)\otimes_{\Q}\CF_{\lambda}$ such that
		\begin{equation}
			\Pi_v(\Gamma_{dR}(A),\Gamma_v(A))=J'_{2,4}\begin{pmatrix}D_v&0\\	0&\frac{1}{2\pi i}(D_v^T)^{-1}\end{pmatrix},
		\end{equation}
		where
		$$	
		J'_{2,4}:=\begin{pmatrix}1&0&0&0\\0&0&0&-1\\0&0&1&0\\0&1&0&0\end{pmatrix},
		$$
		for some $D_v\in\GL_2(\C)$.
	\end{lemma}
	\begin{proof} Set $\CF:=\CF_{\lambda}$ and write 
		$$ 
		\gal(\CF/\Q)=\{\sigma,\tau\},
		$$ 
		where $\sigma=\id_{\CF}$ and $\tau$ is the nontrivial automorphism of $\CF/\Q$. 
		
		Set $V:=H^1_v(A)$ and write $\langle\cdot,\cdot\rangle$ for the Riemann bilinear form on $V$. We use the same notation for its scalar extension to
		$$
		V_{\CF}:=V\otimes_{\Q}\CF.
		$$
		The action of $\CF$ on $V_{\CF}$ gives an eigenspace decomposition 
		\begin{equation}\label{eq:qmsplittingbetti} 
			V_{\CF}=V_{\sigma}\oplus V_{\tau}, 
		\end{equation} 
		where $\dim_{\CF}V_{\chi}=2$ and $\CF$ acts on $V_{\chi}$ via $\chi$, for each $\chi\in\{\sigma,\tau\}$. 
		
		Fix $\gamma,\gamma'\in V_{\sigma}$, and for $e\in\CF$ denote by $e_v$ the endomorphism of $V_{\CF}$ induced by the action of $e$ on $V$. By definition of $V_{\sigma}$, 
		$$
		\langle e_v(\gamma),e_v(\gamma')\rangle = \sigma(e)^2\langle\gamma,\gamma'\rangle. 
		$$ 
		On the other hand, compatibility of the polarization pairing with the Rosati involution gives
		$$
		\langle e_v(\gamma),e_v(\gamma')\rangle=\langle\gamma,(e^{\dagger}\cdot e)\gamma'\rangle,
		$$
		where $\dagger$ denotes the Rosati involution. 
		
		Taking $e=\lambda$ and using $\lambda^\dagger=-\lambda$ from \Cref{qmalglemma1}, comparison of the two identities gives 
		$$
		\langle\gamma,\gamma'\rangle=0. 
		$$ 
		Thus $V_\sigma$ is isotropic. Since $\dim_{\CF}V_\sigma=2$, it follows that $V_\sigma$ is a Lagrangian subspace of $V_{\CF}$. The same argument shows that $V_\tau$ is also Lagrangian.
		
		Together with \eqref{eq:qmsplittingbetti}, this shows that $V_{\sigma}$ and $V_{\tau}$ are transverse Lagrangian subspaces of $V_{\CF}$. We may therefore choose $\gamma_1,\gamma_2\in V_{\sigma}$ and $\delta_1,\delta_2\in V_{\tau}$ such that
		$$
		\Gamma_v(A):=(\gamma_1,\gamma_2,\delta_1,\delta_2)
		$$
		is a symplectic basis of $V_{\CF}$.
		
		Set $V_{dR}:=H^1_{dR}(A/K)$. Since $\CF\subset K$, the $\CF$-action makes $V_{dR}$ an $\CF\otimes_{\Q}K$-module. For $\chi\in\{\sigma,\tau\}$, let $K^\chi$ denote $K$ viewed as an $\CF$-module via $\chi:\CF\rightarrow K$. The decomposition
		$$
		\CF\otimes_{\Q}K\simeq K^\sigma\oplus K^\tau
		$$
		then induces an eigenspace decomposition
		\begin{equation}\label{eq:qmsplitderham}
			V_{dR}= V_{dR}^{\sigma}\oplus V_{dR}^{\tau},
		\end{equation}
		where, for each $\chi\in\{\sigma,\tau\}$, the space $V_{dR}^{\chi}$ is a two-dimensional $K$-vector space on which $\CF$ acts via $\chi$. The same argument as on the Betti side shows that $V_{dR}^{\sigma}$ and $V_{dR}^{\tau}$ are transverse Lagrangian subspaces.
		
		Set	$W:=e_0^{*}\Omega_{A/K}\subset V_{dR},$	where $e_0:\spec(K)\rightarrow A$ denotes the zero section. Since the action of $\CF$ on $V_{dR}$ preserves $W$, the eigenspace decomposition of $V_{dR}$ restricts to
		\begin{equation}\label{eq:qmsplitderham2}
			W=W^{\sigma}\oplus W^{\tau}.
		\end{equation}
		For each $\chi\in\{\sigma,\tau\}$, we have $\dim_K W^\chi=1$, and $\CF$ acts on $W^\chi$ via $\chi$. For the dimension statement, the Skolem--Noether theorem gives an element $u\in D^\times$ whose action interchanges the $\sigma$- and $\tau$-eigenspaces for $\CF$. Since the action of $D$ preserves the Hodge filtration, we have $u(W)=W$. Thus $u$ induces an isomorphism $W^\sigma\simeq W^\tau$, and hence $\dim_K W^\sigma=\dim_K W^\tau=1$.
		
		Choose nonzero vectors $\omega\in W^\sigma$ and $\omega'\in W^\tau$, and extend them to bases 
		$$ 
		(\omega,\eta')\qquad\text{and}\qquad(\omega',\eta)
		$$ 
		of $V_{dR}^\sigma$ and $V_{dR}^\tau$, respectively. Since $W$ is Lagrangian, 
		$$ 
		\langle\omega,\omega'\rangle=0. 
		$$ 
		Moreover, since $V_{dR}^\sigma$ is Lagrangian, we have $\langle\omega,\eta'\rangle=0$. Nondegeneracy of the polarization pairing therefore implies that $\langle\omega,\eta\rangle\neq0$. The same argument gives $\langle\omega',\eta'\rangle\neq0$. 
		
		After replacing $\eta$ by $\eta+c\omega'$ for a suitable $c\in K$, we may assume that
		$$
		\langle\eta,\eta'\rangle=0.
		$$
		Rescaling $\eta$ and $\eta'$ then gives
		$$
		\langle\omega,\eta\rangle=\langle\omega',\eta'\rangle=1.
		$$
		Thus $\Gamma_{dR}(A):=(\omega,\omega',\eta,\eta')$ is a Hodge basis of $H^1_{dR}(A/K)$.
		
		Write $(a_{i,j}):=\Pi_v(\Gamma_{dR}(A),\Gamma_v(A))$. Fix $e\in\CF$. Since $\omega$ is a $\sigma$-eigenvector,
		\begin{equation}
			\rho_v(e^*\omega)=\rho_v(\sigma(e)\omega)=\sigma(e)	(a_{1,1}\gamma_1+a_{1,2}\gamma_2+a_{1,3}\delta_1+a_{1,4}\delta_2).
		\end{equation}
		On the other hand, the $\CF$-equivariance of $\rho_v$ gives
		\begin{multline}
			\rho_v(e^*\omega)=e_v\bigl(a_{1,1}\gamma_1+a_{1,2}\gamma_2+a_{1,3}\delta_1+a_{1,4}\delta_2\bigr)\\	
			=a_{1,1}\sigma(e)\gamma_1+a_{1,2}\sigma(e)\gamma_2+a_{1,3}\tau(e)\delta_1+a_{1,4}\tau(e)\delta_2.
		\end{multline}
		Taking $e=\lambda$ and using $\sigma(\lambda)\neq\tau(\lambda)$, comparison of the two identities gives $a_{1,3}=a_{1,4}=0$.
		
		The same argument for the remaining elements of $\Gamma_{dR}(A)$ gives
		\begin{equation}
			\Pi_v(\Gamma_{dR}(A),\Gamma_v(A))=\begin{pmatrix}
				a_{1,1}&a_{1,2}&0&0\\
				0&0&a_{2,3}&a_{2,4}\\
				0&0&a_{3,3}&a_{3,4}\\
				a_{4,1}&a_{4,2}&0&0
			\end{pmatrix}.
		\end{equation}
		Equivalently,
		\begin{equation}
			\Pi_v(\Gamma_{dR}(A),\Gamma_v(A))=J'_{2,4}	\begin{pmatrix}	D_v&0\\0&D_v'\end{pmatrix}.
		\end{equation}
		
		Write $M:=\Pi_v(\Gamma_{dR}(A),\Gamma_v(A))$. The Riemann bilinear relations give
		\begin{equation}
			M\begin{pmatrix}0&I\\-I&0\end{pmatrix}M^T=\frac{1}{2\pi i}\begin{pmatrix}0&I\\-I&0\end{pmatrix}.
		\end{equation}
		Substituting the block decomposition above gives $D_v'=\frac{1}{2\pi i}(D_v^T)^{-1}$, as required.
	\end{proof}
	\begin{remark}\label{remarkgpcarchimedean}
		Write $\pi_{i,j}$ for the entries of $\Pi_v(\Gamma_{dR}(A),\Gamma_v(A))$. 
		
		Grothendieck's period conjecture predicts that
		$$
		\trdeg_{\bar{\Q}}\bar{\Q}(\pi_{i,j})=4.
		$$
		The above description of $\Pi_v(\Gamma_{dR}(A),\Gamma_v(A))$ gives the upper bound
		$$
		\trdeg_{\bar{\Q}}\bar{\Q}(\pi_{i,j})\leq 5.
		$$
		
		Using a second element $\mu\in D$ with $\mu\lambda\neq\lambda\mu$ yields further relations among the entries of $D_v$ and gives the sharper bound
		$$
		\trdeg_{\bar{\Q}}\bar{\Q}(\pi_{i,j})\leq 4.
		$$
	\end{remark}
	
	\subsubsection{Ordinary places}\label{section:qmperiodsord}
	
	We now assume that $v\in\Sigma_{K,f}$ and first consider the case where $v$ is an ordinary place for $A$.
	\begin{lemma}\label{qmnonarch}
		Let $v\in\Sigma_{K,f}$ be an ordinary place for $A$. Then there exists a symplectic basis $\Gamma_v(A)$ of $H^1_v(A)\otimes_{K_{v,0}}\C_v$ such that, with respect to the $\CF_{\lambda}$-Hodge basis	$\Gamma_{dR}(A)$ chosen in \Cref{qmarch}, we have
		\begin{equation}\label{eq:qmnonarch1}
			\Pi_v(\Gamma_{dR}(A),\Gamma_v(A))=J'_{2,4}\begin{pmatrix}D_v&0\\0&(D_v^T)^{-1}\end{pmatrix}	J'_{2,4},
		\end{equation}
		where $J'_{2,4}$ is as in \Cref{qmarch} and $D_v\in\GL_2(\C_v)$.
	\end{lemma}
	\begin{proof}
		Set $q:=|k_v|$, let $p$ be the residue characteristic, write $q=p^f$, and let $\phi$ denote the $q$-power Frobenius acting on $H^1_v(A)$.
		
		Since $v$ is ordinary for $A$, we have a decomposition
		\begin{equation}\label{eq:frobsplitting}
			H^1_v(A)_{\C_v}:=H^1_v(A)\otimes_{K_{v,0}}\C_v=W_1\oplus W_0,
		\end{equation}
		where, for $j\in\{0,1\}$, the space $W_j$ is spanned by the $\phi$-eigenvectors with eigenvalues $\alpha$ satisfying $v(\alpha)=jv(q)$. Here $\dim_{\C_v}W_j=2$ for each $j\in\{0,1\}$.
		
		Set $\CF:=\CF_\lambda$. Over $\C_v$ we have
		$$
		\CF\otimes_{\Q}\C_v	\simeq	\C_v^{\sigma}\oplus \C_v^{\tau},
		$$
		where $\C_v^{\chi}$ denotes $\C_v$ viewed as a $\CF$-module via the embedding $\chi:\CF\hookrightarrow\C_v$. The induced $\CF$-action therefore gives a decomposition
		\begin{equation}\label{eq:Fsplitcrys}
			H^1_v(A)_{\C_v}=V^\sigma\oplus V^\tau,
		\end{equation}
		where $\CF$ acts on $V^\chi$ via $\chi$.
		
		By the Skolem--Noether theorem, there exists $u\in D^\times$ such that $ueu^{-1}=\tau(e)$ for every $e\in\CF$. The induced action of $u$ on $H^1_v(A)_{\C_v}$ therefore interchanges $V^\sigma$ and $V^\tau$. Since $\dim_{\C_v}H^1_v(A)_{\C_v}=4$, we have
		$$
		\dim_{\C_v}V^\sigma=\dim_{\C_v}V^\tau=2.
		$$
		
		For each $e\in\CF$, functoriality of crystalline cohomology implies that the induced endomorphism $e_{\crys}\in\End(H^1_v(A))$ commutes with $\phi$. Thus, after extending the action of $e_{\crys}$ to $H^1_v(A)_{\C_v}$, the spaces $V^\sigma$ and $V^\tau$ are $\phi$-stable. Hence the $\CF$-eigenspace decomposition is compatible with \eqref{eq:frobsplitting}, giving
		\begin{equation*}
			W_j=W_j^\sigma\oplus W_j^\tau,
		\end{equation*}
		where
		$$
		W_j^\chi:=W_j\cap V^\chi.
		$$
		
		The induced action of $u$ on $H^1_v(A)_{\C_v}$ also commutes with $\phi$, and thus preserves each $W_j$. Since it interchanges the $\sigma$- and $\tau$-eigenspaces, it gives an isomorphism $W_j^\sigma\simeq W_j^\tau$. Hence
		$$
		\dim_{\C_v}W_j^\sigma=\dim_{\C_v}W_j^\tau=1,
		$$
		since $\dim_{\C_v}W_j=2$.
		
		Choose nonzero vectors
		$$
		\gamma\in W_1^\sigma,\qquad\gamma'\in W_1^\tau,\qquad\delta\in W_0^\tau,\qquad\delta'\in W_0^\sigma.
		$$
		Iterating \cite[(2.1.2), p.~118]{dworkogusjacobians} $f$ times, where $q=p^f$, gives the relation
		$$
		\langle\phi(\epsilon),\phi(\epsilon')\rangle=q\langle\epsilon,\epsilon'\rangle.
		$$
		Together with the slope decomposition \eqref{eq:frobsplitting} this shows that both $W_1$ and $W_0$ are isotropic. Since each has dimension $2$, they are Lagrangian, and the pairing between $W_1$ and $W_0$ is perfect.
		
		As in the proof of \Cref{qmarch}, the $\sigma$- and $\tau$-eigenspaces for the $\CF$-action are Lagrangian. Hence
		$$
		\langle\gamma,\delta'\rangle=\langle\gamma',\delta\rangle=0,\qquad\langle\gamma,\delta\rangle\neq0,\qquad\langle\gamma',\delta'\rangle\neq0.
		$$
		After rescaling $\delta$ and $\delta'$, we may assume that
		$$
		\langle\gamma,\delta\rangle=\langle\gamma',\delta'\rangle=1.
		$$
		Thus $\Gamma_v(A):=(\gamma,\gamma',\delta,\delta')$ is a symplectic basis of $H^1_v(A)_{\C_v}$.
		
		As in the proof of \Cref{qmarch}, the $\CF$-equivariance of the comparison isomorphism $\rho_v$, together with the eigenspace decompositions above, forces
		\begin{equation}
			\Pi_v(\Gamma_{dR}(A),\Gamma_v(A))=\begin{pmatrix}a_{1,1}&0&0&a_{1,4}\\0&a_{2,2}&a_{2,3}&0\\	0&a_{3,2}&a_{3,3}&0\\a_{4,1}&0&0&a_{4,4}\end{pmatrix}.
		\end{equation}
		Equivalently,
		\begin{equation}\label{eq: period matrix finite block}
			\Pi_v(\Gamma_{dR}(A),\Gamma_v(A))=J'_{2,4}\begin{pmatrix}D_v&0\\0&D_v'\end{pmatrix}J'_{2,4}
		\end{equation}
		for some $D_v,D_v'\in\GL_2(\C_v)$.
		
		The de Rham--crystalline comparison isomorphism is compatible with the polarization pairings; see \cite[Theorem~5.1.8]{berthbreenmessing}. Hence
		$$
		\Pi_v(\Gamma_{dR}(A),\Gamma_v(A))	\begin{pmatrix}	0&I\\-I&0\end{pmatrix}\Pi_v(\Gamma_{dR}(A),\Gamma_v(A))^T=\begin{pmatrix}0&I\\-I&0\end{pmatrix}.
		$$
		Substituting \eqref{eq: period matrix finite block} gives $D_v'=(D_v^T)^{-1}$, proving the lemma.
	\end{proof}
	
	The preceding lemma can be rewritten as follows.
	\begin{cor}\label{qmnonarchcor}
		In the setting of \Cref{qmnonarch}, there exists $P_v\in\GL_2(\C_v)$ such that, with respect to the bases chosen there,
		\begin{equation}\label{eq:qmnonarch2}
			\Pi_v(\Gamma_{dR}(A),\Gamma_v(A))=J_0^{-1}\begin{pmatrix}	P_v&0\\	0&\dfrac{1}{d_0}P_v	\end{pmatrix}J_0,
		\end{equation}
		where
		$$
		d_0:=\det P_v\qquad\text{and}\qquad J_0:=\begin{pmatrix}1&0&0&0\\0&0&0&1\\0&1&0&0\\0&0&1&0\end{pmatrix}.
		$$
	\end{cor}
	\begin{remark}
		Let $(\pi_{i,j}):=\Pi_v(\Gamma_{dR}(A),\Gamma_v(A))$. 
		
		In contrast to the archimedean case, see \Cref{remarkgpcarchimedean}, in the ordinary case \Cref{qmnonarchcor} gives the sharper bound
		$$
		\trdeg_{\bar{\Q}}\bar{\Q}(\pi_{i,j})\leq 4.
		$$
		
		The action of an element $\mu\in D$ with $\mu\lambda\neq\lambda\mu$ yields the further bound
		$$
		\trdeg_{\bar{\Q}}\bar{\Q}(\pi_{i,j})\leq 3.
		$$
		
		An analogous drop in transcendence degree in the $p$-adic setting also appears for CM elliptic curves; see \cite{andremots}.
	\end{remark}
	\subsubsection{Supersingular places}\label{section:qmsupersingularperiods}
Let $v\in\Sigma_{K,f}$ be a supersingular place for $A$, let $p$ be its residue characteristic, and let $\widetilde{E}_v/k_v$ be the corresponding supersingular elliptic curve from \Cref{reductionqmsurface}. 

The proof of \cite[Theorem~2.2.2]{patankarthesis} also shows that
$$
\End^0_{k_v}(\widetilde{E}_v)=\End^0_{\overline{\mathbb{F}}_p}(\widetilde{E}_v)=D_{p,\infty},
$$
where $D_{p,\infty}$ is the unique quaternion algebra over $\Q$ ramified precisely at $p$ and $\infty$. Hence
$$
\End^0_{k_v}(\widetilde{A}_v)=M_2(D_{p,\infty})=:D_v'.
$$
Let $L(p)$ be an imaginary quadratic field which splits $D_{p,\infty}$.

Set $\CF:=\CF_\lambda$. By our assumption following \Cref{qmlem1}, we have $\CF\subseteq K$. Write $\gal(\CF/\Q)=\{\sigma,\tau\}$, where $\sigma=\id_\CF$ and $\tau$ is the nontrivial automorphism. We denote by $\widetilde{\lambda}_v$ the endomorphism of $H_v^1(A)_{\C_v}$ induced by the reduction $\widetilde{\lambda}_v:\widetilde{A}_v\to\widetilde{A}_v$.

Using \cite[Proposition~5.3.3]{andremsj} together with the $\CF$-action, we obtain the following.
\begin{lemma}\label{qmsupersinglembasis}
	Let $A$, $\lambda$, $v$, $p$, and $L(p)$ be as above. Then there exists a basis $(\gamma_1,\gamma_2,\delta_1,\delta_2)$ of $H_{\crys}^{1}\left(\widetilde{A}_v/W(k_v)\right)\otimes_{W(k_v)}\C_v$ such that:
	\begin{enumerate}
		\item $H_B^1(A_v,\CF L(p)):=\Span_{\CF L(p)}(\gamma_1,\gamma_2,\delta_1,\delta_2)$ is a $4$-dimensional $\CF L(p)$-vector space stable under the action of $D_v'$,
	
		\item $\langle\gamma_1,\gamma_2\rangle=\langle\delta_1,\delta_2\rangle=0$ and
		$$
		\gamma_1\wedge\gamma_2\wedge\delta_1\wedge\delta_2=1;
		$$
		\item $\widetilde{\lambda}_v\gamma_j=\sqrt{\lambda^2}\gamma_j$ and $\widetilde{\lambda}_v\delta_j=-\sqrt{\lambda^2}\delta_j$.
	\end{enumerate}
\end{lemma}
\begin{proof}
	By \cite[Proposition~5.3.3]{andremsj}, for every splitting field $F$ of $D_{p,\infty}$ there exists a $4$-dimensional $F$-vector space
	$$
	H_B^1(A_v,F)\subset H_v^1(A)_{\overline{\Q}_p}
	$$
	stable under $D_v'$. Take $F=\CF L(p)$ and write
	$$
	H_B^1(A_v,\CF L(p))\subset H_v^1(A)_{\overline{\Q}_p}
	$$
	for the resulting space.
	
	Since $\CF\hookrightarrow D_v'$, its action preserves $H_B^1(A_v,\CF L(p))$. As in \Cref{qmnonarch}, the $\CF$-action gives a decomposition
	$$
	H_B^1(A_v,\CF L(p))=V_B^\sigma\oplus V_B^\tau,
	$$
	where $\CF$ acts on $V_B^\chi$ via $\chi$ and $\dim_{\CF L(p)}V_B^\chi=2$ for each $\chi\in\{\sigma,\tau\}$.
	
	Choose bases $(\gamma_1,\gamma_2)$ and $(\delta_1,\delta_2)$ of $V_B^\sigma$ and $V_B^\tau$, respectively. Then condition $(3)$ holds and, as in the proof of \Cref{qmarch},
	$$
	\langle\gamma_1,\gamma_2\rangle=\langle\delta_1,\delta_2\rangle=0.
	$$
	
	By the normalization in the proof of \cite[Proposition~5.3.3]{andremsj}, these bases may moreover be chosen so that
	$$
	\gamma_1\wedge\gamma_2\wedge\delta_1\wedge\delta_2=1.
	$$
\end{proof}

The analogue of \Cref{qmnonarch} in the supersingular case is the following.
\begin{prop}\label{qmsupersingperiod0}
	Let $v$ be a supersingular place for a QM abelian surface $A$ defined over a number field $K$ satisfying the assumptions of \Cref{qmlem1}.
	
	With respect to the $\CF_\lambda$-Hodge basis $\Gamma_{dR}(A)$ chosen in \Cref{qmarch}, we have
	$$
	\Pi_v\left(\Gamma_{dR}(A),\Gamma_v(A)\right)=J_{2,4}'\begin{pmatrix}M_v&0\\0&M_v'\end{pmatrix},
	$$
	where $\Gamma_v(A)$ is the basis chosen in \Cref{qmsupersinglembasis}, $J_{2,4}'$ is as in \Cref{qmarch}, and $\det M_v\det M_v'=1$.
\end{prop}

\begin{proof} 
	Let $(\gamma_1,\gamma_2,\delta_1,\delta_2)$ be as in \Cref{qmsupersinglembasis}. 
	
	Recall the eigenspace decomposition in \eqref{eq:Fsplitcrys},
	$$
	H^1_v(A)_{\C_v}=V^\sigma\oplus V^\tau,
	$$
	where $\mathcal F$ acts on $V^\chi$ via $\chi$.
	By the proof of \Cref{qmsupersinglembasis}, we have $V_B^\chi\subset V^\chi$ for each $\chi\in\{\sigma,\tau\}$. Since $V^\sigma$ and $V^\tau$ are transverse Lagrangians, choose a basis $(v_1,v_2)$ of $V^\tau$ such that $\Gamma_v'(A):=(\gamma_1,\gamma_2,v_1,v_2)$ is symplectic.
	
	As in the proof of \Cref{qmarch},
	$$
	\Pi_v\left(\Gamma_{dR}(A),\Gamma_v'(A)\right)=J_{2,4}'\begin{pmatrix}M_v&0\\0&(M_v^T)^{-1}\end{pmatrix}.
	$$
	By the choice of $\Gamma_v(A)$ in \Cref{qmsupersinglembasis},
	$$
	[\id]_{\Gamma_v'(A)}^{\Gamma_v(A)}=\begin{pmatrix}I_2&0\\0&M'\end{pmatrix}
	$$
	for some $M'\in\GL_2(\C_v)$. Condition $(2)$ of \Cref{qmsupersinglembasis} implies that $\det M'=1$. Hence
	$$
	\Pi_v\left(\Gamma_{dR}(A),\Gamma_v(A)\right)=J_{2,4}'\begin{pmatrix}M_v&0\\0&(M_v^T)^{-1}M'\end{pmatrix}.
	$$
	Set $M_v':=(M_v^T)^{-1}M'$. Then $\det M_v\det M_v'=1$.
\end{proof}

	\section{Background on G-functions}\label{section:backgroundgfuns}

We recall the construction of the $G$-functions used in the sequel and their cohomological interpretation. We refer to \cite[\S~3]{papaspadicpart1} for a more detailed account in the present setting.

\subsection{$G$-functions associated with a family of abelian surfaces}\label{section: multiplce centers}

Let $S$ be a smooth geometrically irreducible curve defined over a number field $K$, let $f:\CX\rightarrow S$ be a one--parameter family of principally polarized abelian surfaces, and let $s_0\in S(K)$. After replacing $f:\CX\rightarrow S$ by its pullback along a finite cover of $S$, and after a finite extension of $K$, we retain the same notation for the resulting data. Following finitely many such replacements, we may assume that there exists $x\in K(S)$ such that the induced morphism $x:S\rightarrow \mathbb{P}^1$ satisfies
\begin{enumerate}
	\item $x(s_0)=0$ and every zero of $x$ is simple and defined over $K$.
	
	\item $x:S\rightarrow \mathbb{P}^1$ is Galois; see \cite[Lemma~5.1]{daworr4}. Moreover, all fibers $\CX_\xi$ with $x(\xi)=0$ are isomorphic to $\CX_{s_0}$.
	
	\item For every $\xi\in S$ with $x(\xi)=0$, the fiber $\CX_\xi$ satisfies the properties of \Cref{qmlem1}.
\end{enumerate}

Fix a zero $\xi$ of $x$ and choose a basis of sections of $H^1_{dR}(\CX/S)$ over an affine open neighborhood $U_{\xi}$ of $\xi$. This gives a matrix of $G$-functions
$Y_\xi\in \GL_4(\bar{\Q}[[x]])$; see \cite{andre1989g}. For a place $v\in \Sigma_K$, let $\Delta_v\subset U_\xi^{\an}$ be an embedded analytic disk, in either the complex analytic or rigid analytic sense, centered at $\xi$ and containing no other zero of $x$.

The basis of sections may be chosen so that the resulting $G$-functions satisfy the following additional property; see \cite[\S~3.3.2]{papaspadicpart1}.
 \begin{center}
	If $|x(s)|_w<\min\{1,R_w(Y_\xi)\}$ for some $s\in\Delta_v\cap S(\bar{\Q})$ and some $w\in\Sigma_{K(s),f}$ lying above $v$, then $\CX_{s}$ has the same reduction modulo $w$ as $\CX_{\xi}$.
\end{center}

The family of $G$-functions $Y_G$ associated with $(\CX\rightarrow S,s_0)$ is obtained by selecting certain matrices $Y_\xi$ as $\xi$ ranges over the zeros of $x$. In more detail, one defines an equivalence relation on the set $\{\xi\in S(K):x(\xi)=0\}$. This equivalence relation records the generic isogenies between families obtained by twisting $\CX\rightarrow S$ by automorphisms of the cover $x:S\rightarrow\mathbb{P}^1$; see \cite[\S~5]{daworr4}. The family $Y_G$ is obtained by choosing one $Y_\xi$ from each equivalence class.

\subsection{Cohomological interpretation of values of $G$-functions}

Let $(f:\CX\rightarrow S,s_0)$ and $x$ be as above. Fix $s\in S(\bar{\Q})$, a zero $\xi$ of $x$, and a place $w\in\Sigma_{K(s)}$ such that $|x(s)|_w<\min\{1,R_w(Y_\xi)\}$. Let $v\in\Sigma_K$ be the place below $w$, and assume that $s\in\Delta_v$. We will say that $s$ is $w$-adically close to $\xi$ if it satisfies these conditions.

The purpose of the next section is to describe relations among the $w$-adic values $\iota_w(Y_\xi(x(s)))$ when the fiber $\CX_s$ is a QM abelian surface. The construction of these relations hinges on the cohomological interpretation of the values $\iota_w(Y_\xi(x(s)))$. This interpretation is provided by Grothendieck's relative de Rham--Betti comparison isomorphism and by the relative de Rham--crystalline comparison isomorphism of Berthelot--Ogus; see \cite{bertogus,ogus}.

Over $\Delta_v$ there is an isomorphism
\begin{equation}
	H^1_{dR}(\CX/S)_{\Delta_v}\rightarrow (H^1_{dR}(\CX/S)_{\Delta_v})^{\nabla}\otimes \CO_{\Delta_v},
\end{equation}
identifying $H^1_{dR}(\CX/S)_{\Delta_v}$ with its space of horizontal sections tensored with $\CO_{\Delta_v}$. Using the corresponding comparison isomorphism over $\Delta_v$, we may identify $(H^1_{dR}(\CX/S)_{\Delta_v})^{\nabla}$ with $H^1_v(\CX_{\xi})$, the Betti or crystalline cohomology of the central fiber; see \cite{daworrpap} and \cite[\S~3]{papaspadicpart1}.

Let $\CP_v(s)$ denote the matrix of the fiber at $s$ of the above isomorphism, with respect to the fixed basis of sections of $H^1_{dR}(\CX/S)$ and a choice of basis of $H^1_v(\CX_\xi)$. Under the above identification, one then has
$$
\CP_v(s)=\iota_w(Y_\xi(x(s)))\cdot \CP_v(\xi).
$$

Assume that the family $\CX\rightarrow S$ is Hodge generic, meaning that the image of the induced morphism $\iota:S\rightarrow\mathcal{A}_2$ is a Hodge generic curve. We may then describe the trivial relations among the elements of $Y_G$, in the terminology of \cite[Chapter~VII]{andre1989g}.

Let $I(\SP_4)$ denote the ideal cutting out the symplectic group $\SP_4$ as a subgroup of $\GL_4$. The ideal of trivial relations is generated by one copy of $I(\SP_4)$ for each of the equivalence classes introduced above; see \cite[Proposition~3.12]{papaspadicpart1}.

In the next section we work in the simplified setting where $x$ has a unique zero, namely $s_0$. In this case, the trivial relations arise by specializing elements of
$I(\SP_4)\subset \bar{\Q}[x,X_{i,j}:1\leq i,j\leq 4]$ at $\zeta:=x(s)$ .

We present both the construction of the relations and the proof of the height bound in this simplified setting. For the relation construction, the passage to the general case is explained in \Cref{section:generalcentres}. For the additional details required in the general case of the height-bound argument, we refer to \cite{papaspadicpart1}.
	
	\section{Relations among values of G-functions}\label{section:qmfiniteplaces}

Throughout this section, let $f:\CX\rightarrow S$ be a one--parameter family of principally polarized abelian surfaces defined over a number field $K$, and fix $s_0\in S(K)$ such that $\CX_0:=\CX_{s_0}$ is a QM abelian surface. After replacing $S$ by an affine open neighborhood of $s_0$, we make the following choices and assumptions:
\begin{enumerate}
	\item $\CX_{0}$ satisfies the properties in \Cref{qmlem1}, 
	
	\item writing $D_0:=\End_{K}^{0}(\CX_{0})$, we fix a presentation $D_0=\left(\frac{a_0,b_0}{\Q}\right)$ with standard basis $\{1,i_0,j_0,i_0j_0\}$, and choose $\lambda_0$ as in \Cref{qmalglemma1}, and
	
	\item the $\CF_{\lambda_0}$-Hodge basis $\Gamma_{dR}(\CX_{0})$ of \Cref{qmarch} extends to a Hodge basis $\Gamma_{dR}(\CX/S)$ of $H_{dR}^{1}(\CX/S)$ consisting of global sections.
\end{enumerate}

Write $\Gamma_{dR}(\CX/S)=(\omega,\omega^{\prime},\eta,\eta^{\prime})$. Its specialization at $s_0$ is the $\CF_{\lambda_0}$-Hodge basis
$\Gamma_{dR}(\CX_0)=(\omega_0,\omega_0^{\prime},\eta_0,\eta_0^{\prime})$ of \Cref{qmarch}.

We retain the notation and conventions of \Cref{section:backgroundgfuns}. In particular, write $x:S\rightarrow\mathbb{P}^{1}$ for the fixed parameter and $Y_G(x)=(y_{i,j}(x))_{1\leq i,j\leq 4}$ for the associated matrix of $G$-functions. For every $s\in S(\bar{\Q})$ and every $v\in\Sigma_{K(s)}$ for which $s$ is $v$-adically close to $s_0$, we have
$$
\Pi_{v}(\Gamma_{dR}(\CX/S)_{s},\Gamma_{v}(\CX_{0}))=\iota_{v}(Y_{G}(x(s)))\, \Pi_{v}(\Gamma_{dR}(\CX_{0}),\Gamma_{v}(\CX_{0})).
$$
At finite places, the choice of $Y_G$ ensures that if $s$ is $v$-adically close to $s_0$, then $\CX_s$ and $\CX_0$ have the same reduction modulo $v$; see \cite[\S~3.3.2]{papaspadicpart1}.

The goal of this section is to construct non-trivial relations with coefficients in $\bar{\Q}$ among the entries of $\iota_v(Y_G(x(s)))$ for points $s\in S(\bar{\Q})$ for which $\CX_s$ is a QM abelian surface. We treat archimedean, supersingular, and ordinary places separately.

\subsection{Relations at archimedean places}\label{section:qmarchrelations}

\begin{prop}\label{qmarchrelations}
	Let $s\in S(\bar{\Q})$ be such that $\CX_s$ is a QM abelian surface. Then there exists a homogeneous polynomial $R_{s,\arch} \in \bar{\Q}\left[X_{i, j}: 1 \leq i, j \leq 4\right]$ such that
	
	\begin{enumerate}
		\item the coefficients of $R_{s,\arch}$ lie in a finite extension $L_s/K(s)$ of degree bounded by an absolute constant,
		\item $\deg R_{s,\arch}\leq c_\infty [K(s):\Q]$, for some positive absolute constant $c_\infty$,
		\item $\iota_{v}(R_{s,\arch}(Y_{G}(x(s))))=0$ for all $v \in \Sigma_{L_{s}, \infty}$ for which $s$ is $v$-adically close to $s_{0}$, and
		\item $R_{s,\arch} \notin I(\SP_{4})$.
	\end{enumerate}
\end{prop}

	\begin{proof} 
By \cite{silverberg}, there exists a finite extension $L_s/K(s)$ of degree bounded by an absolute constant such that $\End_{L_s}^{0}(\CX_s)=\End_{\bar{\Q}}^{0}(\CX_s)$. 

Set $D_s:=\End_{L_s}^{0}(\CX_s)$. Choose a standard basis $\{1,i_s,j_s,i_sj_s\}$ of $D_s$ and set $\alpha_s:=i_s\in D_s\setminus\Q$. Fix $v\in\Sigma_{L_s,\infty}$ such that $s$ is $v$-adically close to $s_0$.

	Functoriality of the de Rham--Betti comparison isomorphism gives
	\begin{equation}\label{eq:alphascommutesarch}
		[\alpha_s]_{dR}\,\Pi_v(\Gamma_{dR}(\CX/S)_s,\Gamma_v(\CX_0))=\Pi_v(\Gamma_{dR}(\CX/S)_s,\Gamma_v(\CX_0))\,[\alpha_s]_v.
	\end{equation}
	Since $\alpha_s$ preserves the Hodge filtration and $\Gamma_{dR}(\CX/S)_s$ is a Hodge basis, we have
	$$
	[\alpha_s]_{dR}=\begin{pmatrix}	A_s&0\\	B_s&C_s	\end{pmatrix},
	$$
	for some $A_s,C_s\in\GL_2(L_s)$ and $B_s\in M_2(L_s)$.
	
Write
$$
\Pi_v(\CX_0):=\Pi_v(\Gamma_{dR}(\CX_0),\Gamma_v(\CX_0)).
$$
Then
$$
\Pi_v(\Gamma_{dR}(\CX/S)_s,\Gamma_v(\CX_0))=\iota_v(Y_G(x(s)))\,\Pi_v(\CX_0).
$$
Since $Y_G$ is generically symplectic, $\det Y_G=1$ and $Y_G^{\adj}=Y_G^{-1}$. Thus \eqref{eq:alphascommutesarch} becomes
	\begin{equation}\label{eq:qmmainarch}
		\iota_v\left(Y_G(x(s))^{\adj}	\begin{pmatrix}	A_s&0\\	B_s&C_s	\end{pmatrix}Y_G(x(s))\right)\Pi_v(\CX_0)=\Pi_v(\CX_0)[\alpha_s]_v.
	\end{equation}

Set
$$
(z_{i,j}(x)):=Y_G(x)^{\adj}\begin{pmatrix}A_s&0\\ B_s&C_s\end{pmatrix}Y_G(x),\qquad Z_{i,j}:=\iota_v(z_{i,j}(x(s))).
$$
By \Cref{qmarch}, \eqref{eq:qmmainarch} becomes
$$
(Z_{i,j})=J_{2,4}^{\prime}\begin{pmatrix}D_v&0\\0&\frac{1}{2\pi i}(D_v^T)^{-1}\end{pmatrix}[\alpha_s]_v\begin{pmatrix}D_v^{-1}&0\\0&(2\pi i)D_v^T\end{pmatrix}(J_{2,4}^{\prime})^{-1}.
$$
Since $(J_{2,4}^{\prime})^{-1}=(J_{2,4}^{\prime})^T$, set
$$
(\widetilde Z_{i,j}):=(J_{2,4}^{\prime})^T(Z_{i,j})J_{2,4}^{\prime}.
$$
Then
$$
(\widetilde Z_{i,j})=\begin{pmatrix}D_v&0\\0&\frac{1}{2\pi i}(D_v^T)^{-1}\end{pmatrix}[\alpha_s]_v\begin{pmatrix}D_v^{-1}&0\\0&(2\pi i)D_v^T\end{pmatrix}.
$$
	
	Write $[\alpha_s]_v=(A_{i,j})_{1\leq i,j\leq 2}$, with $A_{i,j}\in M_2(K)$, by the choice of $\Gamma_v(\CX_0)$ in \Cref{qmarch}. Set
	$$
	(J_{2,4}^{\prime})^T(z_{i,j}(x))J_{2,4}^{\prime}=\begin{pmatrix}\widetilde z_1(x)&\widetilde z_2(x)\\ \widetilde z_3(x)&\widetilde z_4(x)\end{pmatrix},
	$$
	where $\widetilde z_i(x)\in M_2(L_s[[x]])$. Then
	$$
	\iota_v(\widetilde z_2(x(s)))=D_vA_{1,2}(2\pi i)D_v^T,\qquad \iota_v(\widetilde z_3(x(s)))=((2\pi i)D_v^T)^{-1}A_{2,1}D_v^{-1}.
	$$
	Multiplying the two identities and taking traces, we obtain
	$$
	\iota_v\!\left(\tr\!\left(\widetilde z_2(x(s))\widetilde z_3(x(s))\right)\right)=\tr\!\left(D_vA_{1,2}A_{2,1}D_v^{-1}\right)=\tr(A_{1,2}A_{2,1}).
	$$
	
	The preceding identity defines a polynomial $R'_{s}\in L_s[X_{i,j}:1\leq i,j\leq4]$ satisfying
	$$
	R'_{s}(Y_G(x))=\tr\!\left(\widetilde z_2(x)\widetilde z_3(x)\right).
	$$
	Set $t_{0,v}:=\tr(A_{1,2}A_{2,1})\in L_s$ and let
	$$
	H(\underline X)=1-X_{1,1}X_{3,3}-X_{2,1}X_{4,3}+X_{1,3}X_{3,1}+X_{2,3}X_{4,1}\in I(\SP_4).
	$$
	Since $H(Y_G(x))=0$, we obtain
	$$
	\iota_v\!\left(R_{s}^{\prime}(Y_G(x(s)))-t_{0,v}(1-H(Y_G(x(s))))^4\right)=0.
	$$

	Set
	$$
	R_{s,v}:=R_{s}^{\prime}(\underline{X})-t_{0,v}(1-H(\underline{X}))^{4}
	$$
	and
	$$
	R_{s,\arch}:=\prod_v R_{s,v},
	$$
	where the product ranges over all $v\in\Sigma_{L_s,\infty}$ for which $s$ is $v$-adically close to $s_0$, with the convention that the empty product is $1$. By construction, $R_{s,\arch}$ satisfies properties $(1)$--$(3)$ of the proposition. It remains to prove that $R_{s,\arch}\notin I(\SP_4)$.
	
	Suppose $R_{s,\arch}\in I(\SP_4)$. Since $I(\SP_4)$ is prime, some factor $R_{s,v}$ lies in $I(\SP_4)$. Since $(1-H)^4-1\in I(\SP_4)$, this is equivalent to $R'_{s}-t_{0,v}\in I(\SP_4)$. The accompanying Mathematica code \cite{mathematicacodeqm} computes the latter polynomial, denoted there by \texttt{Rqmarch}. In the code,
	$$
	[\alpha_s]_{dR}=\begin{pmatrix}(a_{i,j})&0\\(b_{i,j})&(c_{i,j})\end{pmatrix},\qquad 1\leq i,j\leq2.
	$$
	
	The same code evaluates $R'_{s}-t_{0,v}$ at the symplectic matrix $S(l,m,n,p,q,r)$ encoded by \texttt{valuesForXij} and computes its coefficients as a polynomial in $l,m,n,p,q,r$. Since $R'_{s}-t_{0,v}\in I(\SP_4)$, all these coefficients vanish. The resulting conditions imply
	$$
	b_{1,1}=b_{2,2}=0,\qquad a_{1,2}=a_{2,1}=c_{1,2}=c_{2,1}=0,\qquad a_{1,1}=c_{1,1},\qquad a_{2,2}=c_{2,2}.
	$$
	The coefficient of $l n r^2$ further gives $a_{1,1}=a_{2,2}$. Since $\alpha_s^2\in\Q^\times$, the identity $[\alpha_s]_{dR}^2=\alpha_s^2 I_4$ forces $b_{1,2}=b_{2,1}=0$. The same code also computes the remainder of $R'_{s}-t_{0,v}$ upon division by a Gr\"obner basis of $I(\SP_4)$ after imposing the displayed conditions and $b_{1,2}=b_{2,1}=0$; its vanishing independently gives $a_{1,1}=a_{2,2}$. Hence
	$$
	[\alpha_s]_{dR}=a_{1,1}I_4.
	$$
	Since the action of $D_s$ on $H^1_{dR}(\CX_s/L_s)$ is faithful, it follows that $\alpha_s\in Z(D_s)=\Q$, contradicting the choice of $\alpha_s$.
\end{proof}

	\subsection{Relations at supersingular places}

We now consider finite places. Since $\CX_0$ satisfies \Cref{qmlem1}, it has everywhere good reduction, and we identify $H^1_v(\CX_0)$ with the horizontal sections of $H^1_{dR}(\CX/S)$ near $s_0$ via the de Rham--crystalline comparison isomorphism.
\begin{prop}\label{qmsupersingrel}
	Let $s\in S(\bar{\Q})$ be such that $\CX_s$ is a QM abelian surface. Set
	$$
	\Sigma_{\ssing}(s_0):=\{w\in\Sigma_{K,f}:w\text{ is supersingular for }\CX_0\}
	$$
	and
$$
\begin{aligned}
	\Sigma_{\ssing}(s,s_0):=\{\,w\in\Sigma_{\ssing}(s_0):\;&\exists\,v\in\Sigma_{K(s),f}\text{ with }v\mid w\\
	&\text{such that }s\text{ is }v\text{-adically close to }s_0\,\}.
\end{aligned}
$$
	Then there exists a homogeneous polynomial $R_{s,\ssing}\in\bar{\Q}[X_{i,j}:1\leq i,j\leq4]$ such that
	\begin{enumerate}
\item the coefficients of $R_{s,\ssing}$ lie in an extension $\widetilde{L}_s/L_s$, where $L_s/K(s)$ is the extension of \Cref{qmarchrelations} and $[\widetilde{L}_s:L_s]\leq 2$,
		\item $\deg R_{s,\ssing}\leq c_1[K(s):\Q]|\Sigma_{\ssing}(s,s_0)|$, for some absolute constant $c_1>0$,
		\item $\iota_v(R_{s,\ssing}(Y_G(x(s))))=0$ for all $v\in\Sigma_{\widetilde{L}_s,f}$ for which $s$ is $v$-adically close to $s_0$ and $v\mid w$ for some $w\in\Sigma_{\ssing}(s,s_0)$, and
		\item $R_{s,\ssing}\notin I(\SP_4)$.
	\end{enumerate}
\end{prop}
\begin{proof}
	Let $L_s/K(s)$ be the extension of \Cref{qmarchrelations}, so that $\End_{L_s}^{0}(\CX_s)=\End_{\bar{\Q}}^{0}(\CX_s)$. Set $D_s:=\End_{L_s}^0(\CX_s)$ and choose $\alpha_s\in D_s\setminus\Q$ as in the proof of \Cref{qmarchrelations}. Let $q_s$ be the squarefree positive integer whose prime divisors are exactly the rational primes $p(w)$ lying below the places $w\in\Sigma_{\ssing}(s,s_0)$, set $L_s':=\Q(\sqrt{-q_s})$, and set $\widetilde{L}_s:=L_s L_s'$.
	
	Fix $v\in\Sigma_{\widetilde{L}_s,f}$ such that $s$ is $v$-adically close to $s_0$ and $v\mid w$ for some $w\in\Sigma_{\ssing}(s,s_0)$. Let $\Gamma_v(\CX_0)=(\gamma_1,\gamma_2,\delta_1,\delta_2)$ denote the basis $\Gamma_w(\CX_0)$ of \Cref{qmsupersinglembasis} with\footnote{$L_s'$ is a splitting field of $D_{p(w),\infty}$ for all $w\in \Sigma_{\ssing}(s,s_0)$. } $L(p)=L_s'$.
		
	As in the proof of \Cref{qmarchrelations}, write $[\alpha_s]_{dR}=\begin{pmatrix}A_s&0\\B_s&C_s\end{pmatrix}$. Functoriality of the de Rham--crystalline comparison isomorphism gives
	$$
	[\alpha_s]_{dR}\Pi_v(\Gamma_{dR}(\CX/S)_s,\Gamma_v(\CX_0))=\Pi_v(\Gamma_{dR}(\CX/S)_s,\Gamma_v(\CX_0))[\alpha_s]_v,
	$$
	where $[\alpha_s]_v$ denotes the matrix of the induced action on the crystalline side with respect to $\Gamma_v(\CX_0)$. Hence
	\begin{equation}\label{eq:qmsupermain}
		\iota_v\!\left(Y_G(x(s))^{\adj}\begin{pmatrix}A_s&0\\B_s&C_s\end{pmatrix}Y_G(x(s))\right)=\Pi_v(\CX_0)[\alpha_s]_v\Pi_v(\CX_0)^{-1}.
	\end{equation}
	
	By \Cref{qmsupersingperiod0},
	$$
	\Pi_v(\CX_0)=J_{2,4}^{\prime}\begin{pmatrix}M_v&0\\0&M_v^{\prime}\end{pmatrix},
	$$
	where $M_v,M_v^{\prime}\in\GL_2(\C_v)$ and $\det M_v\det M_v^{\prime}=1$.
	
	Since $\CX_s$ and $\CX_0$ have the same reduction at $v$, the reduction of $\alpha_s$ induces an endomorphism of the common special fiber and hence acts on $H_B^1(\CX_{0,v})$. By the defining property of $H_B^1(\CX_{0,v})$ in \Cref{qmsupersinglembasis}, the matrix $[\alpha_s]_v$ with respect to $\Gamma_v(\CX_0)$ has entries in $\CF_{\lambda_0} L_s'$. Thus
	$$
	[\alpha_s]_v=\begin{pmatrix}A_{1,1}&A_{1,2}\\A_{2,1}&A_{2,2}\end{pmatrix},\qquad A_{i,j}\in M_2(\CF_{\lambda_0} L_s').
	$$
	
	With the notation of the proof of \Cref{qmarchrelations}, \eqref{eq:qmsupermain} gives
	$$
	(\widetilde Z_{i,j})=
	\begin{pmatrix}
		M_vA_{1,1}M_v^{-1}&M_vA_{1,2}(M_v^{\prime})^{-1}\\
		M_v^{\prime}A_{2,1}M_v^{-1}&M_v^{\prime}A_{2,2}(M_v^{\prime})^{-1}
	\end{pmatrix}.
	$$
	It follows that
	$$
	\widetilde Z_2\widetilde Z_3=M_vA_{1,2}A_{2,1}M_v^{-1}.
	$$
	Hence, taking traces,
	$$
	\iota_v(\tr(\widetilde z_2(x(s))\widetilde z_3(x(s))))=\tr(A_{1,2}A_{2,1}).
	$$
	
	Let $R'_s\in \widetilde{L}_s[X_{i,j}:1\leq i,j\leq4]$ be the polynomial defined by
	$$
	R'_s(Y_G(x))=\tr\!\left(\widetilde z_2(x)\widetilde z_3(x)\right),
	$$
	and set $t_{0,v}:=\tr(A_{1,2}A_{2,1})\in \widetilde{L}_s$. Then
	$$
	\iota_v\!\left(R'_s(Y_G(x(s)))-t_{0,v}(1-H(Y_G(x(s))))^4\right)=0,
	$$
	for the same $H\in I(\SP_4)$ as in the proof of \Cref{qmarchrelations}.
	
	Set
	$$
	R_{s,\ssing}:=\prod_v\left(R'_{s}(\underline X)-t_{0,v}(1-H(\underline X))^4\right),
	$$
	where $v$ ranges over the places in $\Sigma_{\widetilde L_s,f}$ for which $s$ is $v$-adically close to $s_0$ and $v\mid w$ for some $w\in\Sigma_{\ssing}(s,s_0)$, with the empty product understood to be $1$. 
	
	Properties $(1)$ and $(3)$ follow directly from the construction. Each factor in the above product has degree at most $8$, while the number of places appearing in the product is at most
	$$
	[\widetilde L_s:\Q]\cdot |\Sigma_{\ssing}(s,s_0)|\le (2[L_s : K(s)] )\cdot [K(s):\Q]\cdot |\Sigma_{\ssing}(s,s_0)|.
	$$
	This proves $(2)$. It remains to prove that $R_{s,\ssing}\notin I(\SP_4)$.
	
	Suppose that $R_{s,\ssing}\in I(\SP_4)$. As in the proof of \Cref{qmarchrelations}, it follows that
	$$
	R'_s-t_{0,v}\in I(\SP_4)
	$$
	for some $v$ as above. This contradicts the non-triviality argument in the proof of \Cref{qmarchrelations}.
	\end{proof}

	\subsection{Relations at ordinary places}\label{section:qmordred}

We now consider finite places of ordinary reduction for $\CX_0$. In contrast to the archimedean and supersingular cases, the polynomials constructed below are independent of the place.

Set $\Sigma_{K,\ord}(s_0):=\{v\in\Sigma_{K,f}:v\text{ is ordinary for }\CX_0\}$.

\begin{prop}\label{propqmfin}
	Let $s\in S(\bar{\Q})$ be such that $\CX_s$ is a QM abelian surface. Then there exists $R_{s,\ord}\in\bar{\Q}[X_{i,j}:1\leq i,j\leq4]$ such that:
	\begin{enumerate}
		\item the coefficients of $R_{s,\ord}$ lie in the extension $L_s/K(s)$ of \Cref{qmarchrelations},
		
		\item $R_{s,\ord}$ is homogeneous of degree at most $8$,
		
		\item $\iota_v(R_{s,\ord}(Y_G(x(s))))=0$ for all $v\in\Sigma_{L_s,f}$ for which $s$ is $v$-adically close to $s_0$ and $v\mid w$ for some $w\in\Sigma_{K,\ord}(s_0)$, and
		
		\item $R_{s,\ord}\notin I(\SP_4)$.
	\end{enumerate}
\end{prop}
\begin{proof} 
Set $D_s:=\End_{L_s}^{0}(\CX_s)$. Fix a presentation
$$
D_s=\left(\frac{a_s,b_s}{\Q}\right)
$$
with standard basis $\{1,\alpha_s,\beta_s,\alpha_s\beta_s\}$, where $\alpha_s^2=a_s$ and $\beta_s^2=b_s$.

Fix $v\in\Sigma_{L_s,f}$ satisfying the conditions in $(3)$, and let $q:=|k_v|$. Since $s$ is $v$-adically close to $s_0$, the fibers $\CX_s$ and $\CX_0$ have the same reduction at $v$. Hence the reduction of $\alpha_s$ acts on $H^1_v(X_0)_{\C_v}$ and, by functoriality of crystalline cohomology, commutes with the $q$-power Frobenius.

Consider the decomposition
$$
H^1_v(\CX_0)_{\C_v}=W_1\oplus W_0
$$
into Frobenius eigenspaces from the proof of \Cref{qmnonarch}. Since $\alpha_s$ commutes with the $q$-power Frobenius, we have $\alpha_s(W_j)\subset W_j$. Hence, with respect to the fixed basis $\Gamma_v(\CX_0)$,
\begin{equation}\label{eq:asmatrix}
	[\alpha_s]_v=
	\begin{pmatrix}
		M_1(\alpha_s)&0\\
		0&M_2(\alpha_s)
	\end{pmatrix},
\end{equation}
for some $M_1(\alpha_s):=(d_{i,j})$ and $M_2(\alpha_s):=(e_{i,j})$ in $\GL_2(\C_v)$.

As in the proof of \Cref{qmsupersingrel}, we obtain \eqref{eq:qmsupermain}. Using \eqref{eq:asmatrix}, the notation of \Cref{qmarchrelations}, and the description of $\Pi_v(\CX_0)$ in \Cref{qmnonarchcor}, this gives
\begin{equation}\label{eq:qmordmain}
	(J_0(Z^{\alpha_s}_{i,j})J_0^{-1})=\begin{pmatrix}	\Pi_0&0\\0&\frac{1}{d_0}\Pi_0\end{pmatrix}J_0\begin{pmatrix}M_1(\alpha_s)&0\\0&M_2(\alpha_s)
	\end{pmatrix}	J_0^{-1}\begin{pmatrix}	\Pi_0^{-1}&0\\0&d_0\Pi_0\end{pmatrix},
\end{equation}
where $\Pi_0=(\pi_{i,j})\in\GL_2(\C_v)$, $d_0=\det\Pi_0$, and $Z^{\alpha_s}_{i,j}:=\iota_v(z^{\alpha_s}_{i,j}(x(s)))$, with $(z^{\alpha_s}_{i,j}(x)):=Y_G(x)^{\adj}[\alpha_s]_{dR}Y_G(x)$.

Set $(\widetilde Z^{\alpha_s}_{i,j}):=J_0(Z^{\alpha_s}_{i,j})J_0^{-1}$. Expanding the right-hand side of \eqref{eq:qmordmain}, we obtain
\begin{equation}\label{eq:qmordrelmain}
	(\widetilde{Z}_{i,j}^{\alpha_s})=\begin{pmatrix}(\pi_{i,j})&0\\0&\frac{1}{d_{0}}(\pi_{i,j})\end{pmatrix}\begin{pmatrix}	d_{1,1} & 0 & d_{1,2} & 0 \\0 & e_{2,2} & 0 & e_{2,1} \\d_{2,1} & 0 & d_{2,2} & 0 \\	0 & e_{1,2} & 0 & e_{1,1}	\end{pmatrix}  \begin{pmatrix}(\pi_{i,j})^{-1} & 0 \\ 0 & d_{0}(\pi_{i,j})^{-1}\end{pmatrix}.
\end{equation}

Write $(\widetilde{Z}_{i,j}^{\alpha_s})=\begin{pmatrix}Q_1^{\alpha_s}&Q_2^{\alpha_s}\\Q_3^{\alpha_s}&Q_4^{\alpha_s}\end{pmatrix}$, where $Q_j^{\alpha_s}\in M_2(\C_v)$. For $\xi,\zeta\in\C_v$, set
$$
\Delta(\xi,\zeta):=(\pi_{i,j})\begin{pmatrix}\xi&0\\0&\zeta\end{pmatrix}(\pi_{i,j})^{-1}.
$$
Since $d_0=\det(\pi_{i,j})$, a direct calculation gives
$$
\Delta(\xi,\zeta)=\begin{pmatrix}\xi+\pi_{1,2}\pi_{2,1}\frac{\xi-\zeta}{d_0}&-\pi_{1,1}\pi_{1,2}\frac{\xi-\zeta}{d_0}\\ \pi_{2,1}\pi_{2,2}\frac{\xi-\zeta}{d_0}&\zeta-\pi_{1,2}\pi_{2,1}\frac{\xi-\zeta}{d_0}\end{pmatrix}.
$$

\begin{claim}
	With notation as above,
	\begin{equation}\label{eq:qmordrel1}
		\widetilde{Z}_{1,2}^{\alpha_s}\widetilde{Z}_{2,3}^{\alpha_s}-\widetilde{Z}_{2,1}^{\alpha_s}\widetilde{Z}_{1,4}^{\alpha_s}=0.
	\end{equation}
\end{claim}
\begin{proof}[Proof of the Claim]
	By \eqref{eq:qmordrelmain}, $Q_1^{\alpha_s}=\Delta(d_{1,1},e_{2,2})$ and $Q_2^{\alpha_s}=d_0\Delta(d_{1,2},e_{2,1})$. We distinguish two cases.
	
	\textbf{Case 1:} $\pi_{1,1}\pi_{1,2}=0$ or $\pi_{2,1}\pi_{2,2}=0$.
	
	By the formula for $\Delta(\xi,\zeta)$, either $\widetilde{Z}^{\alpha_s}_{1,2}=\widetilde{Z}^{\alpha_s}_{1,4}=0$ or $\widetilde{Z}^{\alpha_s}_{2,3}=\widetilde{Z}^{\alpha_s}_{2,1}=0$. Hence \eqref{eq:qmordrel1} holds.
	
	\textbf{Case 2:} $\pi_{i,j}\neq0$ for all $i,j$.
	
	Suppose first that $d_{1,1}=e_{2,2}$. Then $\widetilde{Z}^{\alpha_s}_{1,2}=\widetilde{Z}^{\alpha_s}_{2,1}=0$, so \eqref{eq:qmordrel1} holds. Similarly, if $d_{1,2}=e_{2,1}$, then \eqref{eq:qmordrel1} holds. We may therefore assume that $d_{1,1}\neq e_{2,2}$ and $d_{1,2}\neq e_{2,1}$. By the formula for $\Delta(\xi,\zeta)$,
	$$
	\frac{\widetilde{Z}^{\alpha_s}_{2,1}}{\widetilde{Z}^{\alpha_s}_{1,2}}=\frac{\widetilde{Z}^{\alpha_s}_{2,3}}{\widetilde{Z}^{\alpha_s}_{1,4}}=-\frac{\pi_{2,1}\pi_{2,2}}{\pi_{1,1}\pi_{1,2}},
	$$
	and \eqref{eq:qmordrel1} follows.
\end{proof}

Define $\widetilde{Z}^{\beta_s}_{i,j}$ analogously to $\widetilde{Z}^{\alpha_s}_{i,j}$, with $\beta_s$ in place of $\alpha_s$. Applying the same argument to $\beta_s$ yields
\begin{equation}\label{eq:qmordrel2}
	\widetilde{Z}_{1,2}^{\beta_s}\widetilde{Z}_{2,3}^{\beta_s}-\widetilde{Z}_{2,1}^{\beta_s}\widetilde{Z}_{1,4}^{\beta_s}=0.
\end{equation}
The formula for $\Delta$, applied to the $(1,2)$- and $(2,1)$-entries corresponding to $\alpha_s$ and $\beta_s$, also gives
\begin{equation}\label{eq:qmordrel3}
	\widetilde{Z}_{1,2}^{\alpha_s}\widetilde{Z}_{2,1}^{\beta_s}-\widetilde{Z}_{2,1}^{\alpha_s}\widetilde{Z}_{1,2}^{\beta_s}=0.
\end{equation}

Let $R_{\alpha}$, $R_{\beta}$, and $R_{\alpha\beta}\in L_s[X_{i,j}:1\leq i,j\leq4]$ be the polynomials obtained by expressing the left-hand sides of \eqref{eq:qmordrel1}, \eqref{eq:qmordrel2}, and \eqref{eq:qmordrel3}, respectively, in terms of the entries of $Y_G$. Since $[\alpha_s]_{dR}$ and $[\beta_s]_{dR}$ are independent of $v$, so are $R_{\alpha}$, $R_{\beta}$, and $R_{\alpha\beta}$. We need only show that at least one of these polynomials does not lie in $I(\SP_4)$. Suppose, for contradiction, that all three lie in $I(\SP_4)$.

The accompanying Mathematica code \cite{mathematicacodeqm} computes $R_{\alpha}$, denoted by \texttt{Rqmord0}, and $R_{\alpha\beta}$, denoted by \texttt{Rqmord}. Since $R_{\alpha}\in I(\SP_4)$, the Gr\"obner-basis reduction recorded there gives
$$
\begin{aligned}
	b_{1,1}=b_{2,2}=0,\qquad b_{1,2}=-b_{2,1},\qquad a_{1,1}=c_{1,1},\\
	a_{2,1}=c_{1,2},\qquad a_{1,2}=c_{2,1},\qquad a_{2,2}=c_{2,2}.
\end{aligned}
$$
Similarly, $R_{\beta}\in I(\SP_4)$ yields the same relations for the corresponding uppercase variables $A_{i,j}$, $B_{i,j}$, $C_{i,j}$ in the code, which denote the entries of $[\beta_s]_{dR}$.

After imposing these relations on the variables $a_{i,j}$, $b_{i,j}$, $c_{i,j}$ and $A_{i,j}$, $B_{i,j}$, $C_{i,j}$, the same code gives the coefficients and monomials of the remainder of $R_{\alpha\beta}$ upon division by a Gr\"obner basis of $I(\SP_4)$. The vanishing of these coefficients would imply
$$
[\alpha_s]_{dR}[\beta_s]_{dR}=[\beta_s]_{dR}[\alpha_s]_{dR},
$$
contradicting the relation $\alpha_s\beta_s=-\beta_s\alpha_s$.

If there is no ordinary place $v$ satisfying the conditions in $(3)$, set $R_{s,\ord}:=1$. Otherwise, let $R_{s,\ord}$ be one of $R_{\alpha}$, $R_{\beta}$, and $R_{\alpha\beta}$ which does not lie in $I(\SP_4)$.
\end{proof}

\subsection{The general case}\label{section:generalcentres}

We briefly explain how the preceding discussion extends to the general setting of \Cref{proprelationsintro}. In the construction of \Cref{section:backgroundgfuns}, after passing to the finite cover used to define the family of $G$-functions, the point $s_0$ may have several preimages. Each such preimage $\xi$ serves as a possible center for a matrix of G-functions $Y_{\xi}$ as in \Cref{section: multiplce centers}.

Let $s\in S(\bar{\Q})$ and let $v$ be a place of proximity of $s$ to $s_0$. The construction of G-functions in \cite{papaspadicpart1} is arranged so that a lift of $s$ to the above cover lies in a $v$-adic analytic disk centered at a preimage $\xi$ of $s_0$. The cohomological interpretation of the corresponding matrix of $G$-functions is then exactly the one used in the single-center setting of this section.
	
	\section{Height bounds and supersingular proximity}\label{section:qmina2}

\subsection{The height bound}\label{section: proofofhtbound}	
	\begin{proof}[Proof of \Cref{Theoremmainqmhtbound}]
	The proof is parallel to that of \cite[Proposition~5.1]{papaspadicpart1}. We give the argument in the simplified setting described in \Cref{section:qmfiniteplaces}.
	
	Let $s\in S(\bar{\Q})$ be such that $\CX_s$ is QM. Let $\widetilde{L}_s/K(s)$ be the finite extension introduced in the proof of \Cref{qmsupersingrel}, and set
	$$
	\Sigma(s):=\{v\in \Sigma_{\widetilde L_s}:|x(s)|_v<r_v(Y_G)\},\qquad r_v(Y_G):=\min\{1,R_v(Y_G)\}.
	$$
	Thus $\Sigma(s)$ is the finite set of places of proximity of $s$ to $s_0$. Set
	$$
	\Sigma(s)_{\arch}:=\{v\in\Sigma(s):v\mid\infty\},
	$$
	and define $\Sigma(s)_{\ord}$ and $\Sigma(s)_{\ssing}$ to be the subsets of $\Sigma(s)$ consisting of places at which $\CX_0$ is ordinary and supersingular, respectively. Then
	$$
	\Sigma(s)=\Sigma(s)_{\arch}\sqcup\Sigma(s)_{\ord}\sqcup\Sigma(s)_{\ssing}.
	$$
	If $\Sigma(s)=\emptyset$, then $h(x(s))$, and hence $h(s)$, is bounded by a constant depending only on the fixed data; see the proof of the main height bound in \cite{papasbigboi} for details. We may therefore assume that $\Sigma(s)\neq\emptyset$.
	
	For $\ast\in\{\arch,\ssing,\ord\}$, set $R_{s,\ast}=1$ if $\Sigma(s)_{\ast}=\emptyset$; otherwise, let $R_{s,\ast}$ be the polynomial constructed in \Cref{qmarchrelations}, \Cref{qmsupersingrel}, and \Cref{propqmfin}, respectively. Set
	$$
	R_s:=R_{s,\arch}\cdot R_{s,\ssing}\cdot R_{s,\ord}.
	$$
	
	By \Cref{qmarchrelations}, \Cref{qmsupersingrel}, and \Cref{propqmfin}, the polynomial $R_s$ defines a global relation among the values of the G-functions $Y_G=(y_{i,j}(x))$ at $s$. For every $\ast\in\{\arch,\ssing,\ord\}$, we have $R_{s,\ast}\notin I(\SP_4)$. Since $I(\SP_4)$ is prime, $R_s\notin I(\SP_4)$. Hence the relation is non--trivial in the terminology of \cite[Chapter~VII, \S~5]{andre1989g}. The Andr\'e--Bombieri Hasse principle for values of G-functions \cite[Chapter~VII, \S~5]{andre1989g} therefore gives
	$$
	h(s)\leq c_1\cdot(\deg R_s)^{c_2}.
	$$
	
	By \Cref{qmsupersingrel},
	$$
	\deg(R_{s,\ssing})\leq c_1[K(s):\Q]\cdot|\Sigma_{\ssing}(s,s_0)|.
	$$
	Together with the degree bounds for $R_{s,\arch}$ and $R_{s,\ord}$ from \Cref{qmarchrelations} and \Cref{propqmfin}, this gives
	$$
	\deg R_s\leq c_4[\Q(s):\Q]\left(1+|\Sigma_{\ssing}(s,s_0)|\right)
	$$
	for some constant $c_4>0$. The claimed height bound follows.
\end{proof}

\subsection{Applications to Zilber--Pink}\label{section: zilberpinkapplications}

The Zilber--Pink conjecture in $\mathcal{A}_2$ predicts the following.

\begin{conj}[Zilber--Pink in $\mathcal{A}_2$]\label{zpconjintro}
	Let $C\subset \mathcal{A}_2$ be an irreducible curve defined over $\bar{\Q}$ and not contained in a proper special subvariety. Then
	$$
	\bigcup_{\substack{Y\subset \mathcal{A}_2\ \mathrm{special}\\ \dim Y=1}} C\cap Y
	$$
	is finite.
\end{conj}

As noted in the introduction of \cite{daworr}, there are three types of special curves in $\mathcal{A}_2$:
\begin{enumerate}
	\item $E\times CM$-curves: parameterizing abelian surfaces isogenous to $E\times E'$, where exactly one of the elliptic curves $E$ and $E'$ has complex multiplication;
	
	\item $E^2$-curves: parameterizing abelian surfaces isogenous to $E^2$, for some elliptic curve $E$; and
	
	\item $QM$-curves: parameterizing abelian surfaces with quaternionic multiplication.
\end{enumerate}

A substantial case of \Cref{zpconjintro} was established by Daw and Orr in a series of papers; see \cite{daworr,daworr2,daworr3,daworr5}. Their results assume that the closure of $C$ in the Baily--Borel compactification of $\mathcal{A}_2$ meets the zero-dimensional boundary stratum. From the point of view of the $G$-function method, this provides a natural choice for the distinguished point $s_0$. 

We consider the case where the distinguished point $s_0$ instead corresponds to an intersection of $C$ with a QM special curve. In particular, $s_0$ lies in $\mathcal{A}_2$ itself, rather than on the boundary of its Baily--Borel compactification.

Let $C\subset\mathcal{A}_2$ be a smooth irreducible curve defined over $\bar{\Q}$. Starting by passing to a finite cover carrying the pullback of the universal family, and then making the further finite pullbacks required in the construction of the associated $G$-functions and extending the ground field, we obtain a one--parameter family
$$
f:\CX\longrightarrow S
$$
as in \Cref{section:backgroundgfuns}, together with a finite morphism $\pi:S\rightarrow C$ and a point above $s_0$, which we denote by $\tilde{s}_0$. Write $x$ for the local parameter at $\tilde{s}_0$ chosen as in \Cref{section:backgroundgfuns}.

Let $K$ be a fixed number field of definition of $C$, with $s_0\in C(K)$, over which all the data above are defined. Given $s\in C(\bar{\Q})$, we define
$$
\Sigma_{\ssing}(s,s_0):=\left\{	w\in\Sigma_{K,f} :\exists\,\tilde{s}\in\pi^{-1}(s)\text{ such that }w\in\Sigma_{\ssing}(\tilde{s},\tilde{s}_0)\right\},
$$
where $\Sigma_{\ssing}(\tilde{s},\tilde{s}_0)$ is defined as in	\Cref{Theoremmainqmhtbound}.

Combining \Cref{Theoremmainqmhtbound} with the results of Daw and Orr \cite[\S 6.7]{daworr2} gives the following finiteness result toward the Zilber--Pink conjecture.
\begin{cor}\label{corollaryweakZPQM}
	Let $C\subset\mathcal{A}_2$ be a smooth irreducible curve defined over $\bar{\Q}$ that is not contained in a proper special subvariety. Assume that $C$ contains a QM point $s_0$. 
	
	Then for each $N\in \N$ the set 
	$$
	\{s\in C(\bar{\Q}): s\text{ is QM and }|\Sigma_{\ssing}(s,s_0)|\leq N \}
	$$
	is finite.
\end{cor}

\begin{remark}
	The preceding corollary may be viewed as a Zilber--Pink statement for QM points with controlled integrality. More precisely, let $x$ be the local parameter at $\tilde{s}_0$. Then \Cref{corollaryweakZPQM} gives finiteness for QM points $s$ for which $x(\tilde{s})^{-1}$ is non-integral at places lying above only boundedly many supersingular places of $K$, as $\tilde{s}$ ranges over $\pi^{-1}(s)$. No restriction is imposed at ordinary places.
\end{remark}

\subsection{Places of supersingular proximity}\label{section:placesofssingprox}

The height bound in \Cref{Theoremmainqmhtbound} raises the question of whether the cardinality of the set $\Sigma_{\ssing}(s,s_0)$ of supersingular places of $K$ lying below a place of proximity of $s$ to $s_0$ can be controlled. In \cite{papaspadicpart2}, a generalization by Lauter--Viray \cite{lauterviray} of the work of Gross--Zagier \cite{grosszagier} was used to obtain a purely number--theoretic description of the corresponding set for CM points in one--parameter families of elliptic curves.

For QM abelian surfaces, the set of ordinary places has density $1$ by \cite{sawinordprimes}. By contrast, in the setting of CM elliptic curves studied in \cite{papaspadicpart2}, Deuring's reduction theorem \cite[Chapter~13, Theorem~12]{langellfuns} identifies the ordinary primes with the primes that split in the CM field; hence they have density $\frac{1}{2}$ by Chebotarev.

The QM case of Zilber--Pink in the setting of $\mathcal{A}_2$ would follow from a bound of the form
\begin{equation}\label{eq:boundqmssing}
	|\Sigma_{\ssing}(s,s_0)|\leq c_1(\epsilon)[\Q(s):\Q]^{c_2(\epsilon)}\cdot \disc(\End_{\bar{\Q}}(\CX_s))^{\epsilon},
\end{equation}
for every $\epsilon>0$. Indeed, substituting \eqref{eq:boundqmssing} and the height bound of \Cref{Theoremmainqmhtbound} into the proof of the LGO hypothesis in \cite{daworr2}, and noting that for $\epsilon>0$ sufficiently small the factor $\disc(\End_{\bar{\Q}}(\CX_s))^{\epsilon}$ can be absorbed using the Masser--W\"ustholz upper bound in \cite{mwendoesti}, gives the required estimate.
	
	\subsubsection*{Acknowledgments}
	The author thanks Chris Daw for his encouragement and for many enlightening discussions on G-functions. The author also thanks Or Shahar for showing him the basics on Wolfram Mathematica. The author also heartily thanks the Hebrew University of Jerusalem, the Weizmann Institute of Science, and the IAS for the excellent working conditions.
	
	Work on this project started when the author was supported by Michael Temkin's ERC Consolidator Grant 770922 - BirNonArchGeom. Throughout the majority of this work, the author received funding by the European Union (ERC, SharpOS, 101087910), and by the Israel Science Foundation (grant No. 2067/23). In the final stages of this work, the author was supported by the Minerva Research Foundation Member Fund while in residence at the Institute for Advanced Study for the academic year $2025-26$. 
	
	\subsubsection*{Declaration of AI usage}
	While revising this manuscript from an earlier version available on arXiv, the author used ChatGPT (OpenAI) for editorial and language assistance, principally to improve the exposition. The mathematical ideas, results, proofs, and arguments are solely due to the author.
	
	\bibliographystyle{alpha}
	\bibliography{biblio}
	
\end{document}